\documentclass[12pt,reqno,draft]{amsart}

\usepackage{amscd,amsmath}
\usepackage{amssymb}
\usepackage{amsthm}
\usepackage{enumerate}

\theoremstyle{plain}
\newtheorem{theorem}{Theorem}[section]
\newtheorem{lemma}[theorem]{Lemma}
\newtheorem{proposition}[theorem]{Proposition}
\newtheorem{corollary}[theorem]{Corollary}
\theoremstyle{definition}
\newtheorem{remark}[theorem]{Remark}

\newtheorem{definition}[theorem]{Definition}

\def\zb{\overline{Z}_i}
\def\db{\overline{\partial}}
\def\S{S\mathfrak g^*}
\def\C{\mathbb C}
\def\R{\mathbb{R}}

\def\bu{\bullet}

\renewcommand{\epsilon}{\varepsilon}
\renewcommand{\phi}{\varphi}
\DeclareMathSymbol{\varnothing}{\mathord}{AMSb}{"3F}
\renewcommand{\emptyset}{\varnothing}

\title{Formality in an Equivariant Setting}

\author{Steven Lillywhite}

\address{Department of Mathematics, University of Toronto,
 100 St. George St., Toronto, Ontario,
M5S 3G3} \email{sml@math.toronto.edu}


\begin{document}

\begin{abstract}We define and discuss $G$-formality for certain spaces endowed
with an action by a compact Lie group.  This concept is essentially
formality of the Borel construction of the space in a category of commutative
differential graded algebras over $R=H^\bu(BG)$.  
These results may be applied in
computing the equivariant cohomology of their loop spaces.  
\end{abstract}

\maketitle

\section{Introduction}

In this paper we consider $G$-spaces and give formality results for them
in an equivariant category.  
More specifically, given a $G$-space
$M$, we discuss formality of the Borel construction $EG\times_G M$, or 
equivalently, 
 formality of the complex $A_G^\bu(M)$, of equivariant
differential forms.  
However, in the equivariant setting, the map
$M\to \{pt.\}$ is replaced by $EG\times_G M\to BG$, and consequently
all the commutative
differential graded algebras involved are naturally $R$-algebras, where
$R=H^\bu(BG)$.  Thus formality may be considered in the category of
commutative differential graded $R$ algebras.  We shall also consider
the augmented case, corresponding to equivariant base points which are the
same thing as fixed points of the group action.   
We should like to call a $G$-space $M$ ``equivariantly
formal'' when its Borel construction is formal in the
above sense.  However, 
the term ``equivariant formality'' has come to be used 
to describe the degeneration
of the spectral sequence of the fibration $M\to EG\times_G M\to BG$, owing
to the pervasive influence of \cite{GKM}.  Thus we shall adopt the terminology
``$G$-formal'' in this paper.

We give some general results concerning $G$-formality of products and wedges
and reductions to subgroups.  This is followed by several examples of
$G$-formal spaces, including compact K\"ahler manifolds and formal elliptic
spaces, among others.  Of course, we must make appropriate assumptions on the 
$G$-actions of these spaces for the results to hold.

As an application of these results, we compute the equivariant cohomology
of loop spaces.  (If $M$ is a $G$-space, then so is the loop space of $M$ in
the obvious way).   Our motivation comes
 from considering
the cohomology of symplectic quotients of loop spaces, see \cite{me1},
although the results are of general topological interest.
We shall use an ``equivariant'' bar complex to compute the equivariant
cohomology of the loop space.   If the $G$-space $M$ is $G$-formal,
then the bar complex, which is generally a double complex, loses
a differential and becomes a single complex allowing 
for some easier calculations. In the last section we compute an example.

In an appendix, we discuss bar complexes and Eilenberg-Moore theory
concerning the pull-back of a fibration.   
We also consider equivariant
versions of these results, which are used in several of the proofs
in the main body of the paper.  

In what follows, we shall generally assume that $G$ is a compact, connected
Lie group and that all spaces are connected.  Whenever we have need to use 
the localization theorem in 
equivariant cohomology, we shall assume 
that the spaces under consideration are of the homotopy
type of finite-dimensional $G$-CW-complexes,
and furthermore that they have finitely many connective orbit
types, meaning that the set $\{[G^0_x]\ |\ x\in M\}$ is finite, 
where $G_x$ is  the
stabilizer subgroup at $x$, $G^0_x$ is the connected component of 
the identity, and $[G^0_x]$ denotes the set of conjugacy classes in $G$.
This latter condition is automatically
satisfied, by the way, if $M$ is compact or if $G=S^1$. 

I would like to extend my appreciation to Chris Allday who
took the time to read the manuscript and offered advice on several key
points. In particular, Proposition \ref{GtoT} is due to him.  

\section{$k\mathcal{CDGA}$ and formality}
 
In this section we  recall some important facts about the category
of commutative differential graded algebras,
 the notion of formality, and the connection with
rational homotopy theory.
We shall assume for now 
that our algebras are $k$-algebras, where $k$ is a field of characteristic 
zero.  We shall denote by $k\mathcal{CDGA}^o$
 the category of commutative 
differential graded $k$-algebras which are concentrated in non-negative
degrees and have a differential which raises degree by one.  We assume
further that $H^0(A)\approx k$, for all $A$ in $k\mathcal{CDGA}^o$. We shall
denote by $k\mathcal{CDGA}$ the category of algebras in $k\mathcal{CDGA}^o$
which are augmented over $k$ (i.e. there exists for each $A$ a map
$\epsilon : A\to k$, with $k$ concentrated in degree zero), 
together with augmentation
preserving maps for morphisms. 
We shall call an object of $k\mathcal{CDGA}$ (resp. $k\mathcal{CDGA}^o$)
a $k$CDGA (resp. $k$CDGA$^o$).

We recall Quillen's abstract approach to homotopy theory, \cite{Q1}, 
\cite{Q2}.  
He begins by defining the notion of a closed model category.  A closed model
category is a category, $\mathcal C$,
 with 3 distinguished classes of morphisms called
cofibrations, fibrations, and weak equivalences, which satisfy a number
of axioms.  The homotopy category, $Ho\ \mathcal C$, is defined to be the
localization of $\mathcal C$ with respect to the class of weak equivalences.
Quillen introduces a notion of homotopy and shows that $Ho\ \mathcal C$
is equivalent to the more concrete category $ho\ \mathcal C$ which has
for objects the cofibrant/fibrant objects of $\mathcal C$, and morphisms
the homotopy classes of maps.  We point out the important fact that
two objects $X$ and $Y$ in $Ho\ \mathcal C$
are isomorphic if and only if there exists a chain (in $\mathcal C$)
of weak equivalences

\begin{equation}
X\leftarrow Z_1\to Z_2\leftarrow\dots \leftarrow Z_n\to Y.
\end{equation}

In \cite{BG} it is shown that
the categories $k\mathcal{CDGA}^o$ and $k\mathcal{CDGA}$
are closed model categories where the
weak equivalences are the quasi-isomorphisms (maps which induce an
isomorphism on cohomology),  fibrations are the surjective
morphisms, and cofibrations are maps which satisfy the following 
 lifting condition: a map $f$ is a cofibration if for every commutative
diagram
$$
\begin{CD}
X @>>> V\\
@V{f}VV  @V{p}VV\\
Y @>>> W
\end{CD}
$$
\

\noindent
with $p$ a fibration and weak equivalence, there is a map from $Y$ to
$V$ making the diagram commute. (Actually, in \cite{BG}, the authors
do not assume that $H^0(A)\approx k$ for all algebras $A$.  We have
included this assumption for ease of presentation, but the difference
is slight).

Given a closed model category, $\mathcal C$, with initial object $*$,
an object $B$ is called 
{\it cofibrant} if the map $*\to B$ is a cofibration.
$B$ is called a {\it cofibrant model} for $A$ if $B$ is cofibrant
and there exists a weak equivalence $B\to A$.
It follows from the axioms for a closed model category that
every object in a closed model category has a cofibrant model.  
Moreover, there are various lifting and homotopy results associated with
cofibrant algebras, see \cite{BG} section 6.  We mention one here.  If
$\phi: B_1\to B_2$ is a quasi-isomorphism, and we have a map $f: A\to B_2$ with
$A$ cofibrant, then there exists a lift $\tilde{f} :A\to B_1$ such that 
$\phi\tilde{f}\simeq f$,  where $\simeq$ denotes homotopy.

Note that $k\mathcal{CDGA}$ is pointed with point object $k$. The 
homotopy groups of a $k$CDGA $A$ are defined to be 
$$\pi^nA\overset{def}{=}H^n(\bar A/(\bar A)^2)$$
where $\bar A=\ker\epsilon$, for $\epsilon: A\to k$ a given 
augmentation of $A$.

If $f: B_1\to B_2$ is a weak equivalence of cofibrant $k$CDGA's, then 
$f_*: \pi^\bu B_1\to\pi^\bu B_2$ is an isomorphism.  Thus, if we define
$\Pi^n(A)\overset{def}{=}\pi^n(B)$
 for $B$ a cofibrant model of $A$, then
$\Pi^n(A)$ is well-defined up to isomorphism. Moreover, if $f: A_1\to A_2$
is a map of $k$CDGA's, then $f$ induces a unique homotopy class of maps
$f: B_1\to B_2$, for fixed choices of cofibrant models $B_1, B_2$ of
$A_1, A_2$, respectively.  It follows that there is a unique map
$f_*: \Pi(A_1)\to \Pi(A_2)$.  Thus $\Pi$ is functorial, and
different choices of cofibrant models yield naturally isomorphic
such functors.

In $k\mathcal{CDGA}$, there is a special class of cofibrant models called
minimal models.  A minimal model of an algebra $A$ is defined to be a
cofibrant model, $\mathcal M\to A$, which is connected ($\mathcal M^0\approx
k$), and
such that the induced differential
on $\bar{\mathcal M} /(\bar{\mathcal M})^2$ is zero. 
It can be shown that
each algebra in $k\mathcal{CDGA}$ has a minimal model, unique up to 
isomorphism.  If $M$ is a path-connected topological space, 
the (pseudo-dual) $k$ homotopy groups of $M$ are defined to be: $\Pi^n(M)
\overset{def}{=}\Pi^n(A^\bu(M))=\pi^n(\mathcal M)$, where $\mathcal M$
is a minimal model for $A^\bu(M)$.  Here, 
 $A^\bu(M)$ denotes the Sullivan-de Rham complex which is a $\mathbb Q$CDGA; 
see, for example, \cite{AP} for the definition.   
If $M$ is a smooth manifold, we may also use the
ordinary de Rham complex, taking $k$ to be $\mathbb R$.  
 
Halperin  has explicitly identified the cofibrations (and
hence cofibrant objects) in $k\mathcal{CDGA}$.  Cofibrations are the
so-called KS-extensions, and the cofibrant objects are the KS-complexes.
Since  these notions will be important to us, we give
their definitions here, see \cite{Ha} or \cite{AP}.

\begin{definition}\label{KS}
A map $f: A\to B$ of $k$CDGA's is said to be a {\it KS-extension} if there
exists a well-ordered subset $E\subset B$, $E=\{x_\alpha\}$, 
such that $A\otimes \bigwedge(E)\to B$
is an isomorphism of commutative graded algebras, where $\bigwedge(E)$ denotes
the free graded commutative algebra on the set $E$, and the map is induced
by $f$ and the inclusion of $E\subset B$.  Identifying $B$ with 
$A\otimes \bigwedge(E)$, the differential on $B$ satisfies:
\begin{enumerate}
\item $d_B(a\otimes 1)=d_A(a)\otimes 1$ \label{KS1}
\item $d_B(1\otimes x_\alpha)\in A\otimes \bigwedge(E_{<\alpha})$\label{KS2}
\end{enumerate}
where $E_{<\alpha}=\{x_\beta\ |\ \beta<\alpha\}$.  If $E$ also
satisfies: $deg(x_\alpha)>0 \ \forall x_\alpha\in E$, and $deg(x_\beta)<
deg(x_\alpha)\Rightarrow \beta<\alpha$, then $f$ is called a {\it
minimal KS-extension}.  
If $A=k$,
then we replace the word
``extension'' by the word ``complex'' in the definition, obtaining
the notion of {\it KS-complex}.  (A minimal KS-complex is the same 
thing as a minimal algebra defined above.)
A minimal KS-extension in which $A$ is also
minimal is called a {\it $\Lambda$-minimal $\Lambda$-extension}.  
Note that in a (minimal) KS-extension, $\bigwedge(E)$ is a (minimal) 
KS-complex, with differential such that 
$\epsilon\otimes 1: A\otimes \bigwedge(E)
\to \bigwedge(E)$ is a map of $k$CDGA$^o$'s, where $\epsilon$ is 
the augmentation of $A$.  Moreover, all of these maps may be made 
compatible with augmentations.   

\end{definition}

If $A$ is a $k$CDGA, then its cohomology, $H(A)$, may be considered to
be a $k$CDGA with zero differential.

\begin{definition}$A$ is said to be {\it formal} if $A\approx H(A)$
in $Ho(k\mathcal{CDGA})$.

\end{definition}

\noindent
It is easy to see that this definition is equivalent to the following two.

\begin{lemma}\label{formality}
Consider the category $k\mathcal{CDGA}$. The following are equivalent:
\begin{enumerate}
\item $A$ is formal. \label{formal1}
\item There is a diagram: $$A\leftarrow B\to H(A)$$
where the maps are weak equivalences and $B$ is a cofibrant model for $A$.
(In particular, we may pick $B$ to be minimal).
\label{formal2}
\item There is a chain of quasi-isomorphisms: 
$$A\leftarrow A_1\to A_2\leftarrow \dots \leftarrow A_n\to H(A)$$
\label{formal4}
\end{enumerate}

\end{lemma}

This theory has an important application to rational homotopy theory.  It
turns out that the homotopy category of rational finite $\mathbb Q$-type 
nilpotent spaces is equivalent to the homotopy category of the full 
subcategory of $\mathbb Q\mathcal{CDGA}$ consisting of algebras $A$ with 
$\Pi A$ of finite type, \cite{BG}.
  Thus we may ``do'' rational homotopy theory in a 
category of differential graded algebras.  As an example, if $X$ is
a path-connected, simply-connected topological space of finite 
$\mathbb Q$-type, then there is a natural isomorphism 
$$\Pi^n(A^\bu(X))\approx Hom_\mathbb Q(\pi_n(X)\otimes \mathbb Q, \mathbb Q)
$$ 
where $A^\bu(X)$ is the $\mathbb Q$CDGA of Sullivan-de Rham
differential forms on $X$.  If $X$ is a smooth manifold, the same statement
for homotopy groups 
holds if we use instead the de Rham algebra $A^\bu(X)$ and replace 
$\mathbb Q$ coefficients by $\mathbb R$, or $\mathbb C$.  There is
not a corresponding equivalence of homotopy categories over $\mathbb R$
or $\mathbb C$, however.  

A path-connected topological space is said to be {\it formal} 
if its Sullivan-de
Rham algebra $A^\bu(X)$ is formal.  If $X$ is a smooth manifold, we may
use the de Rham algebra and real or complex coefficients.  However,
a well-known result in rational homotopy theory states that formality over
$\mathbb R$ or $\mathbb C$ implies formality over $\mathbb Q$, see
for example \cite{HalSta}.

Formal spaces
include compact K\"ahler manifolds and many homogeneous spaces including
compact globally symmetric spaces. 
Products, wedges, and connected sums of formal spaces are again formal.  The
topological consequences of formality include the vanishing of
 all Massey products.
 Moreover, the rational homotopy type of such a space is determined
solely by its cohomology algebra (at least for a large class of 
such spaces).

\section{$R\mathcal{CDGA}$ and $G$-formality}
In this paper, we shall be concerned with equivariant versions of standard
formality results.  Let $G$ be a compact, connected Lie group.  Then
$H^\bu(BG; k)$ is isomorphic to the $k$CDGA freely generated by a finite
number of generators of even degree.  We shall denote 
$R\overset{def}{=}H^\bu(BG)$.
We define the category $R\mathcal{CDGA}^o$ to be the category of 
commutative differential graded $R$ algebras.  We shall continue to 
assume that $H^0(A)\approx k$ for all algebras $A$.  
Thus, we obtain a faithful forgetful functor from
$R\mathcal{CDGA}^o$ to  $k\mathcal{CDGA}^o$.   
We also define $R\mathcal{CDGA}$ to be the category of commutative differential
graded $R$ algebras augmented over $R$.  Composing augmentations with 
the augmentation $R\to k$, we get a faithful forgetful functor 
from $R\mathcal{CDGA}$ to $k\mathcal{CDGA}$. 

It is a standard result that if $\mathcal C$ is a closed model category,
$B$ is an object of $\mathcal C$, then the ``over
category'' $\mathcal C/B$ whose
objects are maps $X\to B$ and whose morphisms are commutative squares of
the type
$$
\begin{CD}
X @>{f}>> Y\\
@VVV   @VVV\\
B @= B
\end{CD}
$$

\noindent
may be given the structure of a closed model category with the
following definitions.  Such a morphism
in $\mathcal C/B$ will be called a fibration, cofibration, or weak equivalence,
if the map $f: X\to Y$ is such in $\mathcal C$. 
A similar statement holds for the ``under category'',
$B\backslash C$. See \cite{DS}  for these
and other results about closed model categories. 

Thus we see that both $R\mathcal{CDGA}^o=R\backslash k\mathcal{CDGA}^o$
and $R\mathcal{CDGA}=R\mathcal{CDGA}^o/R$ are closed
model categories. Moreover, the simplicial category structure on 
$k\mathcal{CDGA}^o$ defined in \cite{BG}, section 5, induces a simplicial
category structure on $R\mathcal{CDGA}^o$ and $R\mathcal{CDGA}$ in such a
way that the results of \cite{BG}, section 5, (suitably modified)
hold for these categories as well (cf. \cite{Q1}, II.2, proposition 6).  
From this, it follows that the
homotopy results of \cite{BG}, section 6, (suitably modified) hold
for $R\mathcal{CDGA}^o$ and $R\mathcal{CDGA}$ as well.

\begin{definition} 
We shall say that an $R$CDGA (resp. $R$CDGA$^o$) $A$ is {\it formal}
if  $A\approx H(A)$ in  $Ho (R\mathcal{CDGA})$ (resp. $Ho (R\mathcal
{CDGA}^o$)).

\end{definition} 

If a functor $j: \mathcal C_1\to\mathcal C_2$ between two closed model
categories preserves weak equivalences, then $X\approx Y$ in 
$Ho\ \mathcal C_1$
implies $j(X)\approx j(Y)$ in $Ho\ \mathcal C_2$.  Thus if an algebra $A$
is formal as an $R$CDGA, then it is formal as an $R$CDGA$^o$, and as a 
$k$CDGA, etc.  

Suppose a smooth manifold $M$ has a smooth action of a compact Lie
group $G$.  The equivariant cohomology of $M$ may be computed by means
of the Cartan complex of equivariant differential forms: $A^\bu_G(M)=((S
\mathfrak g^*\otimes A^\bu(M))^G,\ d_G)$ where the differential, $d_G$, 
is zero on
$S\mathfrak g^*$, and for $\alpha \in A^\bu(M)$, $d_G\alpha=d\alpha-
\Sigma u_i\iota_{X_i}\alpha$, where the $\{X_i\}$ are fundamental vector
fields of the action corresponding to a basis of $\mathfrak g$, and the
$\{u_i\}$ are the corresponding algebra generators of $S\mathfrak g^*$,
which are given degree two.
If $M$ is just a topological space, we may compute the equivariant cohomology
of $M$ by means of the $\mathbb Q$CDGA $A_G^\bu(M)$ of \cite{All1}, 
when $G=S^1$.  
Alternatively, we could
use the de Rham algebra of the Borel construction, $A^\bu(EG\times_G M)$
when $M$ is a manifold, or the Sullivan-de Rham algebra of the Borel
construction when $M$ is not a manifold.  We shall let $A_G^\bu(M)$
possibly denote any of the above $k$CDGA's, leaving it to the reader to 
interpret
which model one prefers to use, as well as which ground field $k$.
For a comprehensive treatment of equivariant
de Rham theory, see \cite{GS2}.

Using either model, it is obvious how to obtain an $R$ algebra structure
on $A_G^\bu(M)$.  It is induced by $R\overset{i}{\hookrightarrow}
 A_G^\bu(pt.)\to A_G^\bu(M)$,
where the first map is a choice of minimal model for $A_G^\bu(pt.)$
in $k\mathcal{CDGA}$, and the second map is induced 
from the map $M\to \{pt.\}$.
If we use the Cartan models, then the algebras $A_G^\bu(M)$ are manifestly
augmented over $R$ when the group action has
a fixed point.  This is because in the Cartan model, $A_G^\bu(pt.)= R$,
and the inclusion of a fixed point gives a map $A_G^\bu(M)\to A_G^\bu(pt.)=R$.
However, if we use the Sullivan-de Rham complex of the Borel construction, 
then $A_G^\bu(pt.)=A^\bu(BG)\ne R$.  Thus we must use a quasi-isomorphic
complex which is smaller and augmented over $R$. In \cite{All2}, Allday
shows that the complex $\eta^{-1}(R)$ is quasi-isomorphic to $A_G^\bu(M)$
where $\eta: A_G^\bu(M)\to A_G^\bu(pt.)$ is induced by the inclusion of
a fixed point into $M$, and $R$ is embedded in $A_G^\bu(pt.)$ via $i$ as 
above.   Clearly, $\eta^{-1}(R)$ is augmented over $R$, and is functorial
for equivariant maps of $G$-spaces.  
We shall abuse notation and continue to write $A_G^\bu(M)$, even
when we may really mean $\eta^{-1}(R)$.

Let $\mathcal{GTOP}$ denote the category of path-connected
topological $G$ spaces with morphisms the equivariant maps.  Then the
under category $\{pt.\}\backslash\mathcal{GTOP}$ consists of 
``based G-spaces'',
which is the same thing as $G$ spaces with non-empty fixed-point set
and a choice of basepoint in the fixed-point set.  Then $A_G^\bu(-)$
gives a functor from $\mathcal{GTOP}$ to $R\mathcal{CDGA}^o$
and from $\{pt.\}\backslash\mathcal{GTOP}$ to $R\mathcal{CDGA}$.

\begin{definition}
We shall say that a  $G$-space $M$
is {\it G-formal} if $A_G^\bu(M)$ is formal as an $R$CDGA$^o$.  A
 $G$-space $M$ with equivariant base point $p$ (ie: a choice of fixed point
$p\in M^G$) 
is {\it G-formal at $p$} if $A_G^\bu(M)$ is formal as an $R$CDGA, where
$A_G^\bu(M)$ is augmented via the inclusion of $p$ into $M$. 

\end{definition}

If we continue to define a minimal model of an $R$CDGA  as a
connected 
cofibrant model $\mathcal M$ for which the induced differential on
$\ker\epsilon/(\ker\epsilon)^2$ is zero, where $\epsilon$ is an 
augmentation over $R$,
then there may not
be a minimal model for every algebra in $R\mathcal{CDGA}$.  An
example is $S^1$ acting by rotations of $S^2$ about an axis. It is
easy to see that there can be no minimal model for $A_{S^1}^\bu(S^2)$
in $R\mathcal{CDGA}$. 
However, there is a fairly canonical choice of cofibrant model
for an $R$CDGA.  Let $R\to A$ be an $R$CDGA$^o$.  Then the map $R\to A$,
viewed in $k\mathcal{CDGA}$, may be factored as $R\to R\otimes_k 
\mathcal M\to A$ with the first map the inclusion, the
latter map a quasi-isomorphism, and
$\mathcal M$ a minimal KS-complex, \cite{Ha}. Note that the differential
on $R\otimes\mathcal M$ is not the tensor product differential; see the
definition of a KS-complex above, \ref{KS}. The map $R\to R\otimes\mathcal
M$ is a minimal KS-extension, in particular a cofibration in $k\mathcal{CDGA}$,
and hence we see that $R\otimes\mathcal M$ is a cofibrant model for $A$
in $R\mathcal{CDGA}^o$. 
Suppose $A$ is, moreover, an algebra in $R\mathcal{CDGA}$,
and let $\epsilon: A \to R$ be its augmentation.   
Then composing $R\otimes\mathcal M\to A \overset{\epsilon}{\to} R$ gives an 
$R$-augmentation for $R\otimes\mathcal M$.  
Thus, $R\otimes\mathcal M$ becomes a cofibrant model for $A$ in the
category $R\mathcal{CDGA}$.  As defined, it is unique up to isomorphism.

For those algebras of the form $A_G^\bu(M)$, arising from a group action
on the space $M$, this cofibrant model is more explicitly given by the
Grivel-Halperin-Thomas theorem which states that there is a commutative
diagram

\begin{equation}
\begin{CD}
A_G^\bu(pt.) @>>> A_G^\bu(M) @>>> A^\bu(M) \\
@AAA           @AAA            @AAA     \\
R @>{i}>>         R\otimes_k\mathcal M @>>> \mathcal M 
\end{CD}
\end{equation}

\noindent
associated to the fibration $M\to EG\times_G M\to BG$,
where $\mathcal M$ is a
minimal model for $M$, and the bottom row is a 
$\Lambda$-minimal $\Lambda$-extension, see \cite{Gr}, \cite{Ha}.

\begin{definition}
We shall refer to $R\otimes\mathcal M$ as the {\it $G$-model}
of $A$; or just simply as the {\it $G$-model} of $M$, when $A=A^\bu_G(M)$. 
\end{definition}

\noindent
Sometimes
we may choose to denote it by $\mathcal M_G\overset{def}{=}R\otimes\mathcal M$.
Note that $R\otimes\mathcal M$ may fail to be minimal as a
$k$CDGA.  

Following  \cite{All2}, \cite{AP},
 given a path-connected $G$-space $M$ with equivariant
base-point (ie: a fixed-point) $p$,
the equivariant (pseudo-dual) $k$ homotopy groups are defined 
to be  
\begin{equation} \label{eq homotopy groups}
\Pi^n_{G,p}(M)\overset{def}{=} \pi^n(R\otimes\mathcal M)
= H^n(\ker\epsilon/(\ker\epsilon)^2)
\end{equation}

\noindent
where $\epsilon: R\otimes
\mathcal M \to R$ 
is the  $R$ algebra augmentation  induced by the inclusion of $p$ into $M$,
as above.   
The assignment $(M,p)\mapsto (R\otimes\mathcal
M, \epsilon)$ gives 
a functor from $\{pt.\}\backslash\mathcal
{GTOP}$ to $Ho(R\mathcal{CDGA})$, and  the equivariant 
pseudo-dual $k$-homotopy groups
are functorial as well. Note that if $M$ is $G$-formal, then the equivariant
pseudo-dual $k$ homotopy groups are determined by the equivariant 
cohomology ring of $M$.   

The following lemma is useful for comparing the equivariant 
pseudo-dual $k$ homotopy groups
to the ordinary pseudo-dual $k$ homotopy groups of the Borel construction. 

\begin{lemma}\label{4.2}
Let $A$ be an $R$CDGA and let $R\otimes
\mathcal M$ be the $G$-model for $A$. Then
$R\otimes\mathcal M$ is minimal in $k\mathcal{CDGA}$.  
\end{lemma}
\begin{proof}
We have the augmentation $\epsilon: R\otimes\mathcal M \to R$ which is a map
of $R$CDGA$^o$'s.  The differential, $D$, on $R\otimes\mathcal M$ satisfies
$D(r\otimes 1)=0$, for $r\in R$, and generally has the form
$D(1\otimes\alpha)= r\otimes 1 + \sum r_i\otimes\alpha_i + 1\otimes d\alpha$,
where $\alpha, \alpha_i\in \mathcal M$, $r, r_i\in R$ with $deg(\alpha),
deg(\alpha_i), deg(r), deg(r_i)>0$, and where $d$ is the differential in
$\mathcal M$. Now, $0=\epsilon D(1\otimes\alpha)=
r + \sum r_i\epsilon(\alpha_i) +\epsilon(d\alpha)$. Since $d\alpha\in
(\mathcal M^+)^2$, and $\epsilon$ is an algebra map, it follows that
$\sum r_i\epsilon(\alpha_i) + \epsilon(d\alpha) \in (R^+)^2$.  Hence, we
must have that $r=0$, and it follows that $R\otimes\mathcal M$ is minimal.

\end{proof}

\noindent
As an example, the pseudo-dual $k$ homotopy groups of the 
Borel construction of $S^1$
acting on $S^2$ do not distinguish the trivial action from a standard 
non-trivial one, whereas the equivariant pseudo-dual $k$ homotopy groups do.

\section{Generalities concerning G-formality}

In this section we give some basic results about $G$-formality including
reduction to subgroups and the $G$-formality of products and wedges.

We begin by noting that formality in the category $R\mathcal{CDGA}^o$ 
is equivalent
to formality in $k$CDGA.  In general, for two $R$ algebras $A$ and $B$,
$A\approx B$ in $Ho(k\mathcal{CDGA})$ does not imply that $A\approx B$
in $Ho(R\mathcal{CDGA}^o)$.  Nevertheless, we have the following.

\begin{lemma}
Assume that $R\overset{j}{\to} A$ is an $R$CDGA$^o$ and that we give
$H(A)$ the $R$-algebra structure $R\overset{j^*}{\to}H(A)$.  Then 
$A$ is formal in $k\mathcal{CDGA}$ if and only if $A$ is formal
in $R\mathcal{CDGA}^o$.
\end{lemma}

\begin{proof}
If $A$ is formal in $R\mathcal{CDGA}^o$, then it will be so in 
$k\mathcal{CDGA}$, as we have noted above.  Let us now assume
that $A$ is formal in $k\mathcal{CDGA}$.
Let $\mathcal N$ be a minimal model for $A$ and let  
$R\otimes\mathcal M$ be the $G$-model for $A$.
Then we have a commutative diagram of
$k$CDGA's

\begin{equation}
\begin{CD}
R @>{j}>> A @= A \\
@| @A{\eta}AA  @AAA  \\
R @>>> R\otimes\mathcal M @. \mathcal N\\
@.  @.  @VVV \\
@.  @.  H(A)
\end{CD}
\end{equation}
Since $R\otimes\mathcal M$ is cofibrant in $k\mathcal{CDGA}$, there exists
a map, which is necessarily a quasi-isomorphism,
$R\otimes\mathcal M\to \mathcal N$ making the upper right square
homotopy commute. This gives us a quasi-isomorphism 
$\phi: R\otimes\mathcal M\to \mathcal N \to H(A)$.  Then the map
\begin{equation}
R\otimes\mathcal M\overset{\phi}{\to}H(A)\overset{(\phi^*)^{-1}}{\to}
H(R\otimes\mathcal M)\overset{\eta^*}{\to}H(A)
\end{equation}
is a quasi-isomorphism and a map of $R$-algebras.

\end{proof}

\begin{remark}
We note that this is not true for maps, however.  That is, if
$f: A\to B$ is a map of $R$CDGA$^o$'s, and $f$ is formal as a map
of $k$CDGA's, then $f$ need not be a formal map of $R$CDGA$^o$'s.
\end{remark}

In the category $R\mathcal{CDGA}$, formality is a concept distinct
from formality in the category $k\mathcal{CDGA}$.  In fact, it is easy
to see that $M$ is $G$-formal at $p$ if and only if the map $i: BG\to 
EG\times_GM$ is a formal map, where $i$ is the map induced by the inclusion
of $p$ into $M$.

\begin{definition}
Suppose that $G$ acts on a space $M$.  Then the Serre spectral sequence
associated with the fibration $M\to EG\times_G M\to
BG$ is the same as the spectral sequence (from $E_2$ onwards)
obtained from the $G$-model
$R\otimes\mathcal M$ via the filtration $\mathcal F^p=R^{\ge p}\otimes
\mathcal M$.  If this spectral sequence degenerates at the $E_2$ term,
then \cite{GKM} refers to $M$ as being equivariantly formal.  For obvious
reasons, we wish to avoid this terminology;  however, to
conform as well to current trends, we shall say that $M$ is {\it ef}
when this spectral sequence degenerates at the $E_2$ term.  
\end{definition}

\begin{proposition}\label{sub}
Let $G$ act on a space $M$.  Suppose that $K\subset
G$ is a closed, connected subgroup. 
If $M$ is $G$-formal at $p$ (or $G$-formal) and ef,
then $M$ is $K$-formal at $p$ (resp. $K$-formal).
\end{proposition}

\begin{proof}
We first consider the case where $M$ is $G$-formal at $p$.
The inclusion $K\subset G$ induces a pull-back diagram:

\begin{equation}\label{subsquare}
\begin{CD}
EK\times_K M @>>> EG\times_G M\\
@VVV   @V{p}VV\\
BK @>{i}>> BG
\end{CD}
\end{equation}

We shall denote $H^\bu(BG)$ by $R_G$, and similarly for $R_K$.  
If we are using the Cartan complex of equivariant differential forms,
then there is no problem with the proof.  If we are using Allday's
construction, $\eta^{-1}(R)$, as notated above, then we face the
possibility that this construction may not be functorial with respect 
to changing the group.  This is because  there may not exist  choices
of minimal models so that
$R_G\to R_K\to A^\bu(BK)$ commutes with $R_G\to A^\bu(BG)
\overset{i}{\to}
A^\bu(BK)$.  Then there would
not exist an induced map $\eta^{-1}(R_G)\to \eta^{-1}(R_K)$.

This problem may be circumvented by the following procedure, as pointed out 
to us by C. Allday.  Consider the diagram \ref{subsquare}.  Let 
$\tilde{Y}$ denote the mapping cylinder of the top row, and $Y$ the mapping
cylinder of the bottom row.  Then we have a commutative diagram

\begin{equation}
\begin{CD}
EK\times_K M @>{\tilde{j_1}}>> \tilde{Y} @>{\tilde{j_2}}>> EG\times_G M\\
@VVV  @VVV  @VVV \\
BK @>{j_1}>> Y @>{j_2}>> BG
\end{CD}
\end{equation}

\noindent
in which the maps $j_1, \tilde{j_1}$ induce surjections on differential forms
and the maps $j_2, \tilde{j_2}$ induce quasi-isomorphisms on differential
forms.  It is easy to show that we may use $\tilde{Y}$ to form the complex
$\eta^{-1}(R_G)$, as discussed in section 3, and that this complex will
be quasi-isomorphic to $A^\bu_G(M)$,
and  $G$-formal at $p$ if $M$ is $G$-formal at $p$.  Moreover,
we now may obtain a commutative diagram 

\begin{equation}
\begin{CD}
A^\bu(Y) @>{j_1}>>  A^\bu(BK) \\
@AAA   @AAA  \\
R_G  @>{\tilde i}>>  R_K 
\end{CD}
\end{equation}

\noindent
in which the vertical arrows are quasi-isomorphisms,  
since the map $A^\bu(Y)\to A^\bu(BK)$ is onto.  This follows by the result
for $R\mathcal{CDGA}$ which is the analog of 6.9 in \cite{BG}.

By \ref{applemma1} of the appendix, there is a  quasi-isomorphism
of $k$CDGA's

\begin{equation}
\bar B(A^\bu(BK), A^\bu(BG), A_G^\bu(M))\to A_K^\bu(M)
\end{equation}

\noindent
where we are abusing notation in the event that we are using Allday's
construction. Then, in either case, we obtain a 
quasi-isomorphism of $R_K$CDGA's

\begin{equation}\label{gtok}
\bar B(R_K, R_G, A^\bu_G(M))\to A^\bu_K(M)
\end{equation}

\noindent
The bar complex \ref{gtok} is an $R_K$ algebra via  the
$R_K$ factor, and has an $R_K$-augmentation 
given by $\epsilon(r_K, \alpha)=
r_K\tilde i(\epsilon_G(\alpha))$, where $r_K\in R_K$, $\alpha\in A_G^\bu(M)$,
and $\epsilon_G: A_G^\bu(M)\to R_G$ is the augmentation of $M$ for the
action of $G$.
By the assumption of $G$-formality,
we get a commuting diagram whose vertical arrows are quasi-isomorphisms

\begin{equation}\label{wubba}
\begin{CD}
R_K @<<< R_G  @>>> A_G^\bu(M)\\
@|   @|   @AAA \\
R_K @<<<  R_G @>>>  \mathcal M_G(M) \\
@|  @|  @VVV \\
R_K @<<<  R_G @>>>  H_G^\bu(M)
\end{CD}
\end{equation}

\noindent
Then we obtain the following sequence of maps which are seen to be 
$R_K$CDGA
quasi-isomorphisms by standard comparison theorems for their associated
Eilenberg-Moore spectral sequences.

\begin{equation}
\bar B(R_K, R_G, A_G^\bu(M))\leftarrow
\bar B(R_K, R_G, \mathcal M_G(M))\to
\bar B(R_K, R_G, H_G^\bu(M))
\end{equation}

\noindent
Now the bar complex $\bar B(R_K, R_G, H_G^\bu(M))$ has only the single
differential $\delta$, and computes $Tor_{R_G}(R_K, H_G^\bu(M))$.
Since $M$ is ef, $H_G^\bu(M)$ is a free $R_G$-module.  
Hence we have that $\bar B(R_K, R_G, H_G^\bu(M))_\bu$ is acyclic in
bar degrees greater than zero and the projection to cohomology
\begin{align}
\bar B(R_K, R_G, H_G^\bu(M))_\bu&
\to\bar B(R_K, R_G, H_G^\bu(M))_0\\
&\to H_\bu(\bar B(R_K, R_G, H_G^\bu(M))_\bu) \notag\\
&\approx H_K^\bu(M)\notag
\end{align}
is an $R_K$CDGA quasi-isomorphism.

The case where we consider $M$ to be $G$-formal in the category
$R_G\mathcal{CDGA}^o$ is similar. 
\end{proof}
 
\begin{corollary} \label{little}
Let $G$ act on a
space $M$.  Suppose that $M$ is $G$-formal at $p$ (or $G$-formal)
 and ef.  Then $M$ is formal
in $k\mathcal{CDGA}$.  
\end{corollary}

\begin{proof}
Just take $K$ to be the identity subgroup in \ref{sub}. 

\end{proof}

\begin{remark}\label{Kmodel}
If we use \ref{lastremark} of the appendix, then we can see that a 
$K$-model for $M$ is given by $\bar B_{R_G}(R_K, R_G, \mathcal M_G(M))=
R_K\otimes_{R_G}\mathcal M_G(M)$.
\end{remark}

In line with the general theme of considering maximal tori in compact, 
connected Lie groups, we have the following which is due to C. Allday.
We say that a space $M$ is of {\it finite type} if $H^i(M)$ is a 
finite-dimensional $k$-vector space for all $i$.

\begin{proposition}\label{GtoT}
Let $G$ act on $M$, and let $T\subset G$
be a maximal torus.  If $M$ is $G$-formal ($G$-formal at $p$) then
$M$ is $T$-formal (resp. $T$-formal at $p$).  Moreover, if $M$ is a space of
finite type, $p\in M^G$, and $M$ is $T$-formal at $p$, then $M$ is
$G$-formal at $p$.    
\end{proposition}

\begin{proof}
We can already see that $G$-formal implies $T$-formal by the proof of
\ref{sub}.  We only need the fact that now $R_T$ is a free $R_G$-module
which follows from the well-known fact that as $R_G$-modules,
$R_T\approx R_G\otimes H^\bu(G/T)$.  

Seeing that $T$-formal at $p$
implies $G$-formal at $p$
may be achieved by imitating the  
proof that $A\otimes K$ being formal implies $A$ is formal, for $K$ an 
extension field of $k$, which is
corollary 6.9 of \cite{HalSta}. We omit the details, but mention the setup.
First, we see by \ref{Kmodel} that a $T$-model for $A_T^\bu(M)$ is given by
 $R_T\otimes_{R_G}\mathcal{M}_G(M)$ with differential $1\otimes
D_G$, where $D_G$ is the differential for the $G$-model $\mathcal{M}_G(M)$.
Thus it suffices to show that if $R_T\otimes_{R_G}\mathcal{M}_G(M)$ is
formal as an $R_T$CDGA, then $\mathcal{M}_G(M)$ is formal as an $R_G$CDGA.

It turns out that the constructions of bigraded and filtered models of the 
relevant algebras, as in \cite{HalSta}, give models in the category 
$R\mathcal{CDGA}$.
The proof of corollary 6.9 may be imitated without too much difficulty.

\end{proof}

\begin{proposition}\label{product} 
Suppose that $X$ and $Y$ are  $G$-spaces
which are both $G$-formal (or assume $X$ is $G$-formal at $p$ and $Y$ is
$G$-formal at $q$), 
and suppose that one or both of them
is ef. Then $X\times Y$ is $G$-formal (resp. $G$-formal at $(p,q)$)
for the diagonal action of $G$.
\end{proposition}

\begin{proof}
The pull-back diagram 
\begin{equation}
\begin{CD}
X\times Y @>>> Y \\
@VVV     @VVV  \\
X @>>> \{pt.\} 
\end{CD}
\end{equation}
gives rise to a pull-back diagram
\begin{equation}
\begin{CD}
EG\times_G(X\times Y) @>>> EG\times_G Y \\
@VVV   @VVV \\
EG\times_G X @>>> BG
\end{CD}
\end{equation}

\noindent
Then we obtain an $R$CDGA$^o$ quasi-isomorphism
\begin{equation} 
\bar B(A_G^\bu(X), A_G^\bu(\{pt.\}), A_G^\bu(Y))\overset{\theta}{\to}
A_G^\bu(X\times Y)
\end{equation}

\noindent
by \ref{eqlemma} of the appendix.   If $X$ and $Y$ both
have fixed-points, then so will their product $X\times Y$.  In that
case, $\theta$ is a quasi-isomorphism of $R$CDGA's by \ref{eqlemma}
of the appendix.   Furthermore,

\begin{equation}
\bar B(A_G^\bu(X), R, A_G^\bu(Y))\to 
\bar B(A_G^\bu(X), A_G^\bu(\{pt.\}), A_G^\bu(Y))
\end{equation}

\noindent
is an $R$CDGA$^o$ ($R$CDGA) quasi-isomorphism.  
Since $X$ and $Y$ are $G$-formal, we get $R$CDGA$^o$ ($R$CDGA)
quasi-isomorphisms of bar complexes

\begin{equation}
\begin{CD}
\bar B(A_G^\bu(X), R, A_G^\bu(Y))\\
@AAA\\
\bar B(R\otimes \mathcal M(X), R, R\otimes \mathcal M(Y))\\
@VVV\\
\bar B(H_G^\bu(X), R, H_G^\bu(Y))
\end{CD}
\end{equation}

\noindent
by standard arguments comparing the associated 
Eilenberg-Moore spectral sequences.

Since one or both of $X$, $Y$ is ef,
just as in the proof of \ref{sub},
the bar complex $(\bar B(H_G^\bu(X), R, H_G^\bu(Y))_\bu\ ;\ \delta )$ 
is acyclic in degrees greater than zero with respect to the bar grading,
and the projection to its cohomology
is an $R$CDGA$^o$ ($R$CDGA) quasi-isomorphism.  

\end{proof}

\begin{proposition}\label{wedge}
Let $X$ and $Y$ be  $G$-spaces whose fixed-point sets
are non-emtpy.  Picking base-points  $p\in X^G$ and $q\in Y^G$, 
we may form
the wedge $X\vee Y$ along these base-points. Then $G$ acts on $X\vee Y$.
If $X$ is $G$-formal at $p$ and $Y$ is $G$-formal at $q$, 
then $X\vee Y$ is $G$-formal at the join of $p$ and $q$.
\end{proposition}

\begin{proof}
Let $\epsilon_X, \epsilon_Y$ denote the augmentations of equivariant 
differential forms, and let $i_X, i_Y$ denote the inclusions of $X, Y$
into $X\vee Y$. Then Mayer-Vietoris gives a short exact sequence
\begin{equation}
0\to A_G^\bu(X\vee Y)\overset{i_X+i_Y}{\to} A_G^\bu(X)\oplus A_G^\bu(Y)
\overset{\epsilon_X-\epsilon_Y}{\to} R\to 0.
\end{equation}
\noindent
Thus $i_X+i_Y$ induces an isomorphism $A_G^\bu(X\vee Y)\approx A_G^\bu(X)
\bigoplus_R A_G^\bu(Y)=\ker\{\epsilon_X-\epsilon_Y\}$.  Moreover, since 
$\epsilon_X$, say, induces a surjection
in cohomology, the associated long exact sequence splits into short exact
sequences and thus 
\begin{equation}
H_G^\bu(X\vee Y)\approx H_G^\bu(X)\oplus_R 
H_G^\bu(Y).
\end{equation}
\noindent
Since $X$ and $Y$ are $G$-formal, we have maps
$A_G^\bu(X)\leftarrow \mathcal M_G(X)\to H_G^\bu(X)$ which are
quasi-isomorphisms of $R$CDGA's, and similarly
for $Y$.  
So we have a commutative diagram whose rows are short exact sequences:

\begin{equation}
\begin{CD}
0 @>>> A_G^\bu(X)\oplus_R A_G^\bu(Y) @>>> A_G^\bu(X)\oplus A_G^\bu(Y)
@>>> R @>>> 0 \\
@|  @AAA  @AAA  @|  @| \\
0 @>>> \mathcal M_G(X)\oplus_R\mathcal M_G(Y) @>>> \mathcal M_G(X)\oplus
\mathcal M_G(Y) @>>> R @>>> 0 \\
@| @VVV @VVV @| @| \\
0 @>>> H_G^\bu(X)\oplus_R H_G^\bu(Y) @>>> H_G^\bu(X)\oplus H_G^\bu(Y) 
@>>> R @>>> 0
\end{CD}
\end{equation}

Then we obtain maps between the associated long exact
sequences in cohomology.  By the 5-lemma, it follows that the maps
\begin{equation}
A_G^\bu(X)\oplus_R A_G^\bu(Y)\leftarrow \mathcal M_G(X)\oplus_R
\mathcal M_G(Y) \to H_G^\bu(X)\oplus_R H_G^\bu(Y)
\end{equation}
\noindent  are quasi-isomorphisms.  It is easy to check that these
maps are compatible with augmentations and the $R$ algebra structure,
so are $R$CDGA quasi-isomorphisms.

\end{proof}

\section{Examples of $G$-formal spaces}
In this section we give some examples of $G$-formal spaces.  

\subsection{Compact K\"ahler manifolds}
Let $M$ be a compact
K\"ahler manifold, and $G$ a compact, connected Lie group acting on $M$
 by holomorphic transformations.  
We introduce equivariant holomorphic cohomology groups.  
Since $M$ is a complex manifold, the complex-valued differential forms
on $M$ are bigraded in the usual way. We shall denote $S\mathfrak g^*
\otimes_\R\C$ by simply $S\mathfrak g^*$.  Then we define the equivariant
Dolbeault cohomology to be the cohomology of the complex
\begin{equation}
(\ [S\mathfrak g^*\otimes_{\C}A^{p,\bu}(M)]^G\ ;\ \db+\sum u_i\iota_{Z_i})
\end{equation}
Here $Z_i$ is the holomorphic vector field on $M$ which comes about
by splitting the fundamental vector field $X_i=Z_i+\zb$ into its holomorphic
and anti-holomorphic components.  The generators $u_i\in S\mathfrak g^*$
are given bidegree $(1,1)$. 
The operators act in a similar way
as for the ordinary equivariant cohomology.
We shall denote the qth cohomology of this complex by $H^{p,q}_G(M)$.

The following theorem was proved in \cite{thesis} and independently 
established in \cite{Tele}.

\begin{theorem} \label{main}
Suppose that $M$ is a compact K\"ahler manifold endowed with
a holomorphic action of a compact, connected Lie group $G$, and suppose
that $M$ is ef for the action of $G$.  
Then $M$ is G-formal.  If $M^G\ne\emptyset$, then $M$ is $G$-formal at any
fixed point.  
\end{theorem}
\begin{proof}
The Cartan complex is $(A_G^\bu(M), d_G)=
( [\S\otimes A^{\bu}(M;\C)]^G\ ;\ d+\sum u_i\iota_{X_i})$.  Let $X_i=
Z_i+\zb$ be the splitting of the fundamental vector field $X_i$ into 
its holomorphic and anti-holomorphic parts.  The differential $d=\partial+
\db$ also splits.  Hence we may split the equivariant differential
as $d+\sum u_i\iota_{X_i}=(\db+\sum u_i\iota_{Z_i})+(\partial +\sum u_i
\iota_{\zb})$.
The complex $[\S\otimes A^{\bu}_G(M;\C)]^G$ is bigraded by giving $u_i\in
\S$ bidegree $(1,1)$, and taking the usual bigrading on $A^{\bu}(M;\C)$.
It is easy to show that 
\begin{equation}
([\S\otimes A^{\bu,\bu}(M;\C)]^G\ ;\ (\db+\sum u_i\iota_{Z_i}) , (\partial +
\sum u_i\iota{\zb}))
\end{equation}
is a first quadrant
double complex.
 Accordingly we have two canonical filtrations of this
complex.  We claim that the spectral sequences corresponding to both
of them degenerate at the $E_1$ term, and moreover are n-opposite,
meaning that $'F^p\oplus ''F^q\approx H^n$ for $p+q-1=n$.
Formality for $A^{\bu}_G(M)$ then follows owing to the results in 
\cite{DGMS}, sections 5 and 6.  

Let us consider the filtration where we take $\db+\sum u_i\iota_{Z_i}$
cohomology first.  This is the Dolbeault equivariant cohomology defined
above.  It itself forms a first quadrant
double complex with the two differentials
$\db$ and $\sum u_i\iota_{Z_i}$.  Let us filter so we take $\db$ cohomology
first.  Then the $E_1$ term for the equivariant Dolbeault complex is 
(additively)
\begin{align}
H([\S\otimes A^{\bu}(M)]^G; \db)&\approx [H(\S\otimes A^{\bu}(M); \db)]^G
\approx (\S\otimes H^{\bu}_{\db}(M))^G\\
&\approx (\S)^G\otimes H^{\bu}_{\db}(M)
\approx H^{\bu}(BG)\otimes H^{\bu}_{\db}(M).
\end{align}
 Now by ordinary Hodge theory
for compact K\"ahler manifolds, this last is isomorphic to 
$H^{\bu}(BG)\otimes H^{\bu}(M)$. But now there can be no further non-trivial
differentials in the spectral sequence by the assumption that $M$ is ef.  
This result follows analogously for the other filtration, 
which is just the complex conjugate of this one.  Furthermore, it
is easy to see that the two filtrations are n-opposite. 

Hence we have a ``$\partial_G\db_G$-lemma'' for the equivariant differential
forms, where we mean by $\db_G$ the equivariant Dolbeault operator as defined
above.  Formality follows via the sequence of $\mathbb C$CDGA 
quasi-isomorphisms 
\begin{equation}\label{eqformal}
A_G^\bu(M)\hookleftarrow ker(\db_G) \to H_G^\bu(M)
\end{equation} 
which are the inclusion and projection, respectively.  These maps are maps
of $R$ algebras and moreover, it follows that for equivariant holomorphic
maps between $M$ and $N$, we get a commutative diagram linking  the sequence
\ref{eqformal}
for $M$ to the analogous sequence
for $N$.  In particular, if the action of $G$
on $M$ has fixed-points, then the inclusion of one (chosen as an equivariant
base-point) gives augmentations so that the sequence
\ref{eqformal} commutes with augmentations.  That is, $M$ is $G$-formal in
$R\mathcal{CDGA}$.  
  
\end{proof}

\begin{corollary}
Suppose that $M$ is a compact K\"ahler manifold endowed with a
holomorphic action of a compact, connected Lie group $G$. Assume that
$M^G\ne\emptyset$.  Then $M$ is $G$-formal at any fixed point.
\end{corollary}

\begin{proof} 
Let $p\in M^G$.
Let $T\subset G$ be a maximal torus.  Then $M^T\ne\emptyset$ and a theorem
of Blanchard says that $M$ is ef for the action of $T$, see \cite{Borel},
Chapter XII, theorem 6.2.  By \ref{main}, $M$ is $T$-formal at $p$.  By
\ref{GtoT}, $M$ is $G$-formal at $p$.  
\end{proof}

\begin{remark}
The proof of \ref{main} implies an equivariant Hodge decomposition:
\begin{equation}
H^n_G(M)\approx \bigoplus_{p+q=n}H^{p,q}_G(M)
\end{equation}

\end{remark}

\subsection{Elliptic spaces}
We recall that
an {\it elliptic space} $M$ is a space such that both $H^\bu(M;k)$
and $V$ are finite-dimensional $k$-vector spaces, where $\mathcal M(M)=
\bigwedge(V)$ is a minimal model for $M$.
We shall use the following result of \cite{Lup}.

\begin{proposition}(Lupton)\label{Lupton}
 Let $F\to E\to B$
be a fibration in which $F$ is formal and elliptic, and $B$ is formal
and simply-connected.
If the Serre spectral sequence of the fibration degenerates at the $E_2$ term,
then $E$ is formal also.  

\end{proposition}

\begin{theorem}\label{elliptic}
Let $M$ be an elliptic $G$ space.  If $M$ is formal and ef, then
$M$ is $G$-formal. If $M^G\ne\emptyset$, then $M$ is $G$-formal at any
fixed point. 
\end{theorem}

\begin{proof}
We have the fibration $M\to EG\times_G M\to BG$.  Then \ref{Lupton}
implies that $A_G^\bu(M)$ is formal as a $k$CDGA.  The proof of Lupton's
proposition works (adapting to our situation) by finding a model
for $A_G^\bu(M)$ of the form $R\otimes \mathcal M$ which is bigraded as
a $k$CDGA.  Here, $\mathcal M$ is the bigraded (minimal) model of $M$.
Elements of $R$ are in degree zero for the second grading, so that
$(R\otimes\mathcal M)_0 = R\otimes (\mathcal M)_0$. 
It is shown that with respect to the second grading, 
$H_+(R\otimes \mathcal M)=0$ and hence the projection to cohomology
\begin{equation}
R\otimes\mathcal M\to (R\otimes\mathcal M)_0  \to H_G^\bu(M)
\end{equation}
is a quasi-isomorphism.  Clearly,
this is a map of $R$ algebras.  Moreover, if $M^G\ne \emptyset$, then
this map commutes with the augmentations over $R$.  This follows because
firstly the map $R\otimes\mathcal M\to (R\otimes\mathcal M)_0$ commutes
with augmentations.  Secondly, since
the augmentation $\epsilon: R\otimes\mathcal M \to R$ is map of $R$CDGA's,
$\epsilon(d\alpha)=0$ for all $\alpha$, so that the map 
$(R\otimes\mathcal M)_0\to H_G^\bu(M)$ commutes with augmentations.       

\end{proof} 

\begin{corollary}
Let $M$ be an elliptic space.  Suppose that a torus $T$ acts on $M$
with $M^T\ne\emptyset$.  Suppose further that one of the components
of the fixed-point set, say $M^T_i$, satisfies $H^{odd}(M^T_i)=0$.  
Then $M$ is $T$-formal at any fixed-point.
\end{corollary}

\begin{proof}
Since $M$ is elliptic, it follows (via localization and localization for
equivariant rational homotopy \cite{AP}) that each component of the fixed-point
set is elliptic and $\chi_\pi(M)=\chi_\pi(M^T_i)$ where $\chi_\pi$ is the
homotopy Euler characteristic.  But Halperin has shown that for elliptic
spaces the conditions $H^{odd}=0$ and $\chi_\pi=0$ are equivalent and
moreover such spaces are  formal.  Thus
$0=\chi_\pi(M^T_i)=\chi_\pi(M)$.  Hence $M$ is formal and $H^{odd}(M)=0$.
But this latter condition implies that $M$ is ef.  So we may apply
\ref{elliptic}.

\end{proof}

\begin{remark}
Suppose $G$ acts on a simply-connected space
$M$ with non-empty fixed-point set.  Then by
picking a base-point in the fixed-point set, we obtain an action
of $G$ on the space of based loops in $M$, denoted $\Omega M$. 
Since the cohomology of $\Omega M$ is free, we see that 
$\Omega M$ will be $G$-formal if $\Omega M$ is ef.  (Lupton's
proof could be extended to this case, as well).  If $G=T$ is a torus,
and $M$ is elliptic, then the condition that $\Omega M$ is ef is
equivalent to the $G$-model $R\otimes\mathcal M(M)$ being minimal
in the category $R\mathcal{CDGA}$, see \cite{AP}, 3.3.15.
\end{remark}

\subsection{Miscellania}
Next we shall give a few extra examples of $G$-formality.  

\begin{theorem}\label{even}
Let $M$ be a space with minimal model $\mathcal M=
\bigwedge(V)$.  Suppose that $dx=0$ for all $x\in V^{even}$ such that
$deg(x)< dimM$.  Suppose further that the circle $S^1=T$ acts on $M$, that
$M$ is ef, and that each component of the fixed-point set is formal and
satisfies $H^{odd}(M^T)=0$.  Then
$M$ is $T$-formal at any fixed point.  
\end{theorem}

\begin{proof}
Since $M$ is ef, the Serre spectral sequence for the fibration $M\to 
ET\times_T M \to BT$ degenerates at the $E_2$ term.  (Note that by the
localization theorem, this implies that $M^T\ne\emptyset$). By a standard
change of basis argument, we may assume that in the $T$-model
$(R\otimes\mathcal M, D)$ we have $Dx=0$, for $x\in V^{even}$ such that
$deg(x)< dimM$.  Let $i: M^T \hookrightarrow M$ denote the inclusion of the 
fixed-point set. Then we have  maps of
$R$CDGA's (actually, the algebras on the right-hand side of the diagram do not
satisy $H^0=k$, but this will not present any problems)
\begin{equation}
\begin{CD}
A_T^\bu(M) @>{i}>> A_T^\bu(M^T) \\
@AAA    @AAA   \\
R\otimes\mathcal M(M) @>{i}>> R\otimes\mathcal M(M^T) \\
@.  @V{h}VV  \\
H_T^\bu(M) @>{i^*}>> H_T^\bu(M^T) = R\otimes H^\bu(M^T)
\end{CD}
\end{equation}
\noindent
where $h$ is a quasi-isomorphism since $M^T$ is formal.   
Since $M$ is ef, the map $i^*$ is an injection. 
We claim that $h i(R\otimes\mathcal M(M))\subseteq i^*(H_T^\bu(M))$. 
Since the maps are algebra maps, it suffices to check this
on algebra generators.   
Since $M$ is ef, the localization theorem shows that
$i^*$ is an isomorphism in degrees $\ge dimM$.   Also if $\alpha\in R\otimes
\mathcal M(M)$ has odd degree, then $hi(\alpha)=0$ since
$H_T^{odd}(M^T)=0$ by assumption.  So it suffices to check 
the claim on algebra
generators of $R\otimes\mathcal M(M)$ of even degree less than $dimM$.
Let $\alpha$ be such a generator.  If $\alpha\in R$, then the claim
is obviously true.  If $\alpha\in \mathcal M(M)$, then by assumption
$D\alpha=0$.  Then $hi(\alpha)=[i(\alpha)]=i^*([\alpha])$.  
Thus we have a map
\begin{equation}
R\otimes\mathcal M(M)\overset{hi}{\to} H_T^\bu(M).
\end{equation}
which is a quasi-isomorphism of $R$CDGA's.

\end{proof}

\begin{corollary}
Let $M$ be a  simply-connected space with minimal model
$\bigwedge(V)$.  Suppose that $dx=0$ for all $x\in V^{even}$ such that
$deg(x)< dimM$.  Suppose further that a torus $T$ acts on $M$,
that $M$ is ef, and that each component of
the fixed-point set is formal and satisfies $H^{odd}(M^T)=0$.
Then $M$ is formal in $k\mathcal{CDGA}$.  
\end{corollary}

\begin{proof}
First of all, there is a subcircle $S^1\subset T$ such that $M^{S^1}=
M^T$.  The inclusion of this circle $S^1\hookrightarrow T$ induces
a pullback diagram:

\begin{equation}
\begin{CD}
ES^1\times_{S^1}M @>>> ET\times_T M\\
@VVV   @VVV\\
BS^1 @>>> BT
\end{CD}
\end{equation}

Since the action of $T$ is ef, the Serre spectral sequence for the
fibration on the right degenerates at the $E_2$ term.  But then the
same is true for the pull-back fibration.  Hence the $S^1$ action is
ef as well.  Now the result follows from \ref{even} and \ref{little}.

\end{proof}
 
\begin{corollary}
Let $M^4$ be a  space such that $H^{odd}(M)=0$,
and $dimM=4$.
Suppose that a circle $S^1=T$ acts on $M$.  Then $M$ is $T$-formal at any fixed
point.
\end{corollary}

\begin{proof}
We have that $H^{odd}(M)=0$ so that $M$ is ef.  Then $M^T\ne\emptyset$.
By localization, $H^{odd}(M^T)=0$.  But path-connected
spaces with $H^1=0$ of dimension less than or equal to 4 are formal, 
so each component of $M^T$ is formal.
The result follows by \ref{even}.

\end{proof}

\begin{remark}
A simple example of an $S^1$-space satisfying the conditions of
\ref{even}, but which is not K\"ahler or elliptic is the following.  Let
$S^1$ act on $S^4$ so that the fixed-point set consists of two isolated
points.  Extend this to a  diagonal action of $S^1$
on $S^4\times S^4$.  Then, removing
a neighborhood of a fixed-point, we may form the connected sum 
$S^4\times S^4\#S^4\times S^4$.  This manifold then inherits an $S^1$ action
with 6 isolated fixed-points.  It is not elliptic and not even symplectic
since $H^2=0$.  It is easy to check that it satisfies the conditions of
\ref{even}, so is $S^1$-formal.  (This can also be seen by 
proving that the connected sum (made in an equivariant setting) of
$G$-formal spaces is again $G$-formal, which we have omitted).  
\end{remark}

We conclude this section with two examples which do not involve the condition
of $M$ being ef.  

\begin{lemma}\label{free}
Let $M$ be a simply-connected compact manifold.  Suppose that 
$G$ acts freely on 
$M$ and $dimG\ge dimM-6$.  Then $M$ is $G$-formal.
\end{lemma}

\begin{proof}
Since $G$ acts freely, $M/G$ is a simply-connected manifold of dimension
6 or less.  Hence $M/G$ is formal \cite{Miller}.  So $EG\times_G M$
is formal. 

\end{proof}

\begin{remark}
Suppose, in the situation of \ref{free}, we have that $dimM-6>rank(G)$. 
Let $T\subset G$ be a maximal torus.  Then by \ref{GtoT},
$M/T$ is a 
simply-connected manifold of dimension greater than 6 which is formal.
\end{remark}

\begin{lemma}
Let $M$ be a simply-connected elliptic space.  Suppose that 
$G$ acts almost-freely on $M$
(meaning all isotropy groups are finite), and  
$rank(G)=-\chi_\pi(M)$.  Then $M$ is $G$-formal.
\end{lemma}

\begin{proof}
Since $M$ and $BG$ have finite-dimensional pseudo-dual rational
homotopy, so does
$EG\times_G M$, as may be seen by considering the fibration 
$M\to EG\times_G M\to BG$.  Since $G$ acts almost freely, $H^\bu(
EG\times_G M)$ is finite-dimensional as well.  
Furthermore
\begin{equation}
\chi_\pi(EG\times_G M)=\chi_\pi(M)+\chi_\pi(BG)=-rank(G)+rank(G)=0.
\end{equation}
Thus $EG\times_G M$ is elliptic with $\chi_\pi=0$, so is formal.   
\end{proof} 

\section{An Application}\label{application}
In this section we give an  application of $G$-formality. We will show
that the computation of the equivariant cohomology of loop spaces
simplifies considerably when the space is $G$-formal.  

Let us consider a simply-connected space $M$.
Suppose that $G$ acts on $M$ with non-empty fixed-point set.
Let $p\in M^G$ be a choice of base-point.  Then we get an action of 
$G$ on the loops in $M$ based at $p$, $\Omega(M;p)$, which we shall
often abbreviate as $\Omega M$.  Let $P(M;p)$ be the space of paths
in $M$, based at $p$.  Then we have the fibration

\begin{equation}
\begin{CD}
\Omega M @>>> P(M;p)\\
@VVV     @V{\pi}VV \\
\{p\} @>>> M
\end{CD}
\end{equation}

\noindent
where $\pi$ is the map sending a path $\gamma(t)$ to its value at time 1,
$\gamma(1)$.
Moreover, the $G$-action induces a pull-back diagram of fibrations

\begin{equation}
\begin{CD}
EG\times_G\Omega M @>>> EG\times_G P(M;p)\\
@VVV  @V{\pi}VV\\
BG @>>> EG\times_G M
\end{CD}
\end{equation}

\noindent
Hence there is a quasi-isomorphism of $R$CDGA's

\begin{equation}\theta: \bar B(R, A_G^\bu(M), A_G^\bu(P(M;p))
\to A_G^\bu(\Omega(M))
\end{equation}  

\noindent
by \ref{eqlemma} of the appendix.   
Now the inclusion of $\{p\}$ into $P(M;p)$ followed by $\pi$ is
the inclusion of $\{p\}$ into $M$.  These  maps are equivariant
so induce their analogs on the Borel constructions.  
Hence we get an $R$CDGA quasi-isomorphism

\begin{equation}
\bar B(R, A_G^\bu(M), R)\leftarrow \bar B(R, A_G^\bu(M), A_G^\bu(P(M;p)).
\end{equation}   

\begin{proposition}
Let $G$ act on a simply-connected space
$M$ with non-empty fixed-point set, so that $G$ acts
on $\Omega M$.  Suppose that $M$ is $G$-formal.  Then there is
an isomorphism of $R$ algebras
\begin{equation}
H_G^\bu(\Omega M)\approx Tor_{H_G^\bu(M)}(R,R)
\end{equation}

\end{proposition}

\begin{proof}
We have that $A_G^\bu(\Omega M)$ is quasi-isomorphic to
 $\bar B(R, A_G^\bu(M), R)$ (via a sequence of $R$CDGA quasi-isomorphisms).
The assumption of $G$-formality means we have a commuting diagram
of $R$-algebras
\begin{equation}
\begin{CD}
R @<{\epsilon}<< A_G^\bu(M) @>{\epsilon}>> R\\
@|  @AAA  @|\\
R @<{\epsilon}<< \mathcal M_G(M) @>{\epsilon}>> R\\
@|  @VVV  @|\\
R @<{\epsilon}<< H_G^\bu(M) @>{\epsilon}>> R
\end{CD}
\end{equation}

\noindent
We obtain $R$CDGA quasi-isomorphisms
\begin{equation}
\bar B(R, A_G^\bu(M), R)\leftarrow \bar B(R, \mathcal M_G(M), R)\to
\bar B(R, H_G^\bu(M), R).
\end{equation}

\noindent
This follows by standard comparison theorems for the Eilenberg-Moore
spectral sequences associated to the bar complexes.    
Thus we have
that $\bar B(R, H_G^\bu(M), R)$ is quasi-isomorphic to $A_G^\bu(\Omega M)$
(via a sequence of $R$CDGA quasi-isomorphisms).  
But the cohomology of
$\bar B(R, H_G^\bu(M), R)$ is  $Tor_{H_G^\bu(M)}(R,R)$.  

\end{proof}

\begin{remark}\label{remarkcomplex} 
We can always choose any resolution to compute $Tor$.   But we
note that we may always use the bar resolution, and using \ref{overR}
of the appendix,
we see that when 
$M$ is $G$-formal, $H_G^\bu(\Omega M)$ may be computed via the
(single) complex 

\begin{equation}
(\bar B_R(R, H_G^\bu(M), R); \delta).
\end{equation}

\end{remark}

\begin{remark}  
We could also obtain analogous results for the equivariant cohomology
of the free loop space $LM$. 
\end{remark}

\section{An Example}
In this section we compute an example of the equivariant cohomology
of the based-loop space using the normalized bar complex over $R$ of
\ref{remarkcomplex}.

\subsection{Example : $S^1$ acting on $\Omega S^2$}
The circle $S^1$ acts on the 2-sphere $S^2$ by rotations about an axis, say
the z-axis when $S^2$ is the unit sphere in $\R^3$. This action is holomorphic
and Hamiltonian. Thus by \ref{main}, $S^2$ is $G$-formal, ($G=S^1$). 
It is easy to show that the equivariant
cohomology ring is
\begin{equation}
H^{\bu}_{S^1}(S^2; k)\approx k[x,u]/(x+u)(x-u)
\end{equation}
where the degree of $x$ and $u$ is two, and $R=k[u]$ acts as multiplication
by $u$.  

The fixed-point set, $F$, consists of the north and south poles.  We shall 
write $F=\{N,S\}$.  Let $\Omega S^2$ be loops based at the north pole.  Then
$S^1$ acts on $\Omega S^2$. 
Then the equivariant cohomology of the based loops, $H_{S^1}(\Omega S^2)$,
may be computed as the cohomology of the bar complex
\begin{equation}
(\bar B_{k[u]}(k[u],\frac{k[x,u]}{(x^2-u^2)},k[u])\ ;\  \delta).
\end{equation}
  Let $\omega$ be
the symplectic form on $S^2$.  Then $x$ is represented by the form $\omega-
uf\in A^{\bu}_{S^1}(S^2)$,
and $u$ is represented by the form $u\in A^{\bu}_{S^1}(S^2)$, using
the Cartan complex of equivariant differential forms. 
Here, $f$ is the moment map which sends a point on $S^2\subset
\mathbb R^3$ to its $z$-component.  
Then the 
inclusion of the north pole $\{N\}$ into $S^2$ induces the 
augmentation $H^{\bu}_{S^1}(S^2)
\to H^{\bu}_{S^1}(\{N\})\approx k[u]$ sending $x\mapsto -u$ and $u\mapsto u$.
We omit the details of computing the bar complex, 
but one finds without difficulty
the cohomology generators $(1,\underbrace{x,\dots,x}_{n},1)$
in degree n for n odd, and $(u^{n/2},1)$ in degree n for n even.  Owing
to the shuffle product structure on the bar complex, one sees that as 
an $R$ algebra,
\begin{equation}
H^{\bu}_{S^1}(\Omega S^2)
\approx k[u]\oplus\bigwedge(x_1)\oplus\bigwedge(x_3)\oplus\bigwedge(x_5)
\dots
\end{equation}
where $x_i$ is an indeterminant of degree $i$.  

\begin{remark}
In this example, the normalized bar complex $\bar B_R(R, H_{S^1}^\bu(S^2), R)$
is actually isomorphic to the $k$CDGA
minimal model for $ES^1\times_{S^1}\Omega
S^2$, which is 
\begin{equation}
\mathcal M_{ES^1\times_{S^1}\Omega S^2}=\bigwedge(u,x,y) \ \ \  
(du=0=dx,\ dy=ux)
\end{equation}
where the degrees of $u$ and $y$ are 2, and the degree of $x$ is 1. 
The isomorphism is given by $(1,x,1)\mapsto x$, $(u,1)\mapsto u$,
and $(1,x,x,1)\mapsto y$.
\end{remark}

\begin{remark}
In this example, the space $ES^1\times_{S^1}\Omega S^2$ is not formal, 
implying that $\Omega S^2$ is not G-formal.  Indeed, Massey
products abound. 
\end{remark} 

\appendix 
\section{Bar complexes and Eilenberg-Moore theory}
In this appendix we shall discuss the theory of Eilenberg and Moore
concerning pull-backs of fibrations. We will also consider equivariant
versions of these results.
For references, see \cite{Mc},  \cite{Smith}, or \cite{EMo}.

Let us suppose that we have a fibration $F\to E\overset{p}
{\to} B$ and a map $f: X\to
B$, so that we obtain a pull-back diagram:

\begin{equation}\label{diagram}
\begin{CD}
E_f @>{\tilde f}>> E \\
@V{\tilde p}VV   @V{p}VV \\
X @>{f}>>  B 
\end{CD}
\end{equation}

\noindent
Then the maps $f^*$ and $p^*$  make $A^\bu(X)$ and $A^\bu(E)$
(differential graded) modules over $A^\bu(B)$. Let us assume that $B$
is simply-connected.  Then a theorem of Eilenberg and Moore asserts
that there is an isomorphism 
\begin{equation}
\theta: Tor_{A^\bu(B)}(A^\bu(X), A^\bu(E))\to H^\bu(E_f).
\end{equation} 

\noindent
We may use the bar resolution to obtain a resolution of, say, $A^\bu(X)$
by $A^\bu(B)$-modules.  Since we are considering $A^\bu(-)$
to be the  de Rham 
or Sullivan-de Rham complex, we will use Chen's
normalized bar resolution, see \cite{Ch} or \cite{GJP}.  

More specifically, the bar complex is

\begin{equation}
B_k(A^\bu(X), A^\bu(B), A^\bu(E))=\bigoplus_{i=0}^{\infty}
A^\bu(X)\otimes_k (sA^\bu(B))^{\otimes i}\otimes_k A^\bu(E)
\end{equation}

\noindent
where the tensor products are over the ground field $k$, and $s$ denotes
the suspension functor on graded vector spaces which lowers degree by one.  
Hence the degree of an element $(\alpha,\omega_1,\dots,\omega_k,\beta)$
is:  $deg(\alpha)+\sum_{i=1}^k (deg(\omega_i)-1)+deg(\beta)$,
where $\alpha\in A^\bu(X)$, $\omega_i\in A^\bu(B)$, and $\beta\in 
A^\bu(E)$. Actually,
the bar complex is bigraded.  We introduce the {\it bar degree}, denoted
$B_k(A^\bu(X), A^\bu(B), A^\bu(E))_{\bu}$.  The bar degree of an element
$(\alpha,\omega_1,\dots,\omega_k,\beta)$ is defined to be $-k$.  The
other grading is the normal tensor product grading, the degree of an
element $(\alpha,\omega_1,\dots,\omega_k,\beta)$ being 
$deg(\alpha)+\sum_{i=1}^k deg(\omega_i)+deg(\beta)$.

There are two differentials of total degree +1:

\begin{align}
d(\alpha,\omega_1,\dots,\omega_k,\beta)&=(d\alpha, \omega_1,\dots,\omega_k,
\beta)\\\notag
&+\sum_{i=1}^k(-1)^{\epsilon_{i-1}+1}
(\alpha,\omega_1,\dots,\omega_{i-1},d\omega_i,\omega_{i+1},\dots,\omega_k,\beta
)\\\notag
&+(-1)^{\epsilon_k}(\alpha,\omega_1,\dots,\omega_k,d\beta)\\
-\delta(\alpha,\omega_1,\dots,\omega_k,n)&=(-1)^{\epsilon_0}
(\alpha\omega_1,\omega_2,\dots,\omega_k,\beta)\\\notag
&+\sum_{i=1}^{k-1}(-1)^{\epsilon_i}(\alpha,\omega_1,\dots,
\omega_{i-1},\omega_i\omega_{i+1},\omega_{i+2},\dots,\omega_k,\beta)\\\notag
&+(-1)^{\epsilon_{k-1}+1}(\alpha,\omega_1,\dots,\omega_{k-1},\omega_k\beta)
\end{align}

\noindent
where 
$\epsilon_i=deg\alpha+
deg\omega_1+\dots +deg\omega_i-i$. The differential $\delta$ has
degree $+1$ with respect to the bar grading, while the differential
$d$ has degree $+1$ with respect to the tensor product grading. 
One may verify that  
$d\delta+\delta d=0$, and we put $D\overset{def}{=}
d+\delta$ to be the total differential.
With the given bigrading,
we get a double complex with the two differentials $d$ and $\delta$
which gives rise to the Eilenberg-Moore spectral sequence.

Chen's normalized version of this bar complex is the following.
  If $f\in A^0(B)$, let $S_i(f)$ be the operator on 
$B(A^\bu(X),A^{\bu}(B),A^\bu(E))$ defined by 
\begin{equation}
S_i(f)(\alpha,\omega_1,\dots,
\omega_k,\beta)
=(\alpha,\omega_1,\dots,\omega_{i-1},f,\omega_i,\dots,\omega_k,\beta)
\end{equation}
for
$1\le i\le k+1$.  Let $W$ be the subspace of $B(A^\bu(X),A^{\bu}(B),A^\bu(E))$
generated by the images of $S_i(f)$ and $DS_i(f)-S_i(f)D$.  Then define

\begin{equation}
\bar B(A^\bu(X),A^{\bu}(B),A^\bu(E))\overset{def}{=}
B(A^\bu(X),A^{\bu}(B),A^\bu(E))/W.
\end{equation}

\noindent  
Then $W$ is closed under $D$ and when $H^0(B)=k$ (B is connected), then
$W$ is acyclic so that $\bar B(A^\bu(X), A^\bu(B), A^\bu(E))$ 
is quasi-isomorphic
to $B(A^\bu(X), A^\bu(B), A^\bu(E))$.  Notice that in the normalized
bar complex, there are no elements of negative degree, and 
with our assumption that $B$ is simply-connected, we are
assured convergence of  the associated Eilenberg-Moore spectral sequence.
The map $\theta$ 
mentioned above is
induced by the map 
\begin{equation}
\theta: B(A^\bu(X), A^\bu(B), A^\bu(E))\to A^\bu(E_f)
\end{equation}
which sends all tensor products to zero except for 
$A^\bu(X)\otimes_k A^\bu(E)$, where the map is: $(\alpha, \beta)\mapsto
\tilde p^*\alpha\wedge \tilde f^*\beta$.  Note that $\theta(W)=0$, so that
we get an induced map 
\begin{equation}
\theta: \bar B(A^\bu(X), A^\bu(B), A^\bu(E))\to
A^\bu(E_f).  
\end{equation}

The normalized bar complex may also be augmented.  The augmentation,
$\epsilon$, maps all elements of positive total degree to zero.  The
elements of degree zero have the form $(f,g)$, where $f\in A^0(X)$
and $g\in A^0(E)$.  Then we define $\epsilon(f,g)=\epsilon_X(f)\epsilon_E(g)
=f(x_0)g(e_0)$,
where $x_0$ and $e_0$ are chosen base-points in $X$ and $E$, respectively,
and $\epsilon_X, \epsilon_E$ are the augmentations of $A^\bu(X), 
A^\bu(E)$, respectively.
If we choose base-points so that the pull-back diagram above preserves
all base-points, then $\theta$ is an augmentation preserving map.

The bar complex has a natural coalgebra structure.  Since we are
inputting $k$CDGA's to the bar complex, we also obtain a structure
of $k$CDGA on the bar complex via the shuffle product.  

More specifically, if $(a_1, \dots, a_p)$
and $(b_1, \dots, b_q)$ are two ordered sets, then a {\it shuffle}
$\sigma$ of $(a_1, \dots, a_p)$ with $(b_1, \dots, b_q)$ is a permutation
of the ordered set $(a_1, \dots, a_p, b_1, \dots, b_q)$ which preserves
the order of the $a_i$'s as well as the order of the $b_j$'s.  That is,
we demand that if $i<j$, then $\sigma(a_i)<\sigma(a_j)$ and 
$\sigma(b_i)<\sigma(b_j)$.  

We  obtain a product on $B(A^\bu(X), A^\bu(B), A^\bu(E))$
by first taking the normal tensor product on the $A^\bu(X)\otimes A^\bu(E)$
factors, then taking the tensor product of this product with the 
shuffle product on the $A^\bu(B)^{\otimes i}$ factors.  As usual,
we introduce a sign $(-1)^{deg(\alpha)deg(\beta)}$ whenever $\alpha$
is moved past $\beta$. One checks that this product induces a product
on Chen's normalized complex, $\bar B(A^\bu(X), A^\bu(B), A^\bu(E))$,
as well.  
Thus we arrive at the following lemma, whose  proof is left
to the reader. For more details, see \cite{me3}.

\begin{lemma}\label{applemma1}
Assume that we have the pull-back diagram  \ref{diagram}, 
where $p$ is a fibration and $B$ is simply connected. Then the normalized
bar complex 
$$
\bar B(A^\bu(X), A^\bu(B), A^\bu(E))
$$
 is a $k$CDGA. 
Moreover, 
$$\theta : \bar B(A^\bu(X), A^\bu(B), A^\bu(E))\to A^\bu(E_f)$$
is a quasi-isomorphism of $k$CDGA's.  
\end{lemma}

\begin{remark}
We note that Chen's normalization is functorial.  That is, if we have
a commutative diagram of $k$CDGA's
\begin{equation}
\begin{CD}
A_2 @<<< B_2 @>>> C_2 \\
@AAA   @AAA   @AAA  \\
A_1 @<<< B_1 @>>> C_1
\end{CD}
\end{equation}
then we get a map of $k$CDGA's $\bar B(A_1, B_1, C_1)\to 
\bar B(A_2, B_2, C_2)$.
\end{remark}

We may formulate an equivariant analog of the bar complex.
Let us consider again the pull-back diagram \ref{diagram}.
If we suppose further that the spaces $X, B$, and $E$ are $G$-spaces, 
and that $f$ and $p$ are equivariant maps, then we obtain a pull-back
diagram:

\begin{equation}\label{eqdiagram}
\begin{CD}
EG\times_G E_f @>{\tilde f}>> EG\times_G E\\
@V{\tilde p}VV    @V{p}VV\\
EG\times_G X @>{f}>> EG\times_G B
\end{CD}
\end{equation}

\noindent 
Note that we are assuming
$B$ to be simply-connected which in turn implies that
$EG\times_G B$ is simply-connected as well. We may apply \ref{applemma1}
to the diagram \ref{eqdiagram}. However, the bar complex
$\bar B(A_G^\bu(X), A_G^\bu(B), A_G^\bu(E))$ has the extra structure
of an $R$CDGA$^o$ or $R$CDGA, depending on fixed points.  
We may give it an $R$ algebra strucure via the $R$
algebra structure on the $A_G^\bu(X)$ factor, and we define the augmentation
as above, assuming that we can choose our base-points as described before
to be actually fixed points of the group action. We arrive at the following.

\begin{lemma}\label{eqlemma}
Assume that in the pull-back diagram \ref{diagram},
 we have that $X$, $B$, and $E$ are all
$G$-spaces with $f$ and $p$ equivariant maps. 
Then the normalized
bar complex

$$
\bar B(A_G^\bu(X), A_G^\bu(B), A_G^\bu(E))
$$
is an $R$CDGA$^o$.  Moreover,

$$
\theta: \bar B(A_G^\bu(X), A_G^\bu(B), A_G^\bu(E))\to A_G^\bu(E_f) 
$$
is a quasi-isomorphism of $R$CDGA$^o$'s.  
If we  assume further that all fixed-point
sets are non-empty, and the diagram \ref{diagram} preserves base-points
chosen from the various fixed-point sets, then the normalized bar complex
is an $R$CDGA, and $\theta$ is a quasi-isomorphism of $R$CDGA's.

\end{lemma}

In this equivariant
case, we may further simplify the bar complex, following an idea of \cite{GJP}.
Let us consider the bar complex over $R$:
\begin{equation}
B_R(A_G^\bu(X), A_G^\bu(B), A_G^\bu(E))=\bigoplus_{i=0}^{\infty}
A_G^\bu(X)\otimes_R (sA_G^\bu(B))^{\otimes i}\otimes_R A_G^\bu(E)
\end{equation}
where all the tensor products are over $R$.  

\begin{lemma}\label{overR}
Suppose that $A, B$, and $C$ are $R$CDGA's and we have morphisms of
$R$CDGA's $A\leftarrow B\to C$, where $R=H^\bu(BG)$ for $G$ a compact,
connected Lie group.  (We use this sequence to define a (differential graded)
$B$-module structure on $A$ and $C$). Suppose further 
either that for each $r\in R$,
$r$ is not a zero-divisor in $A$, or that this condition holds for
$C$.
Then the natural projection

\begin{equation}
B_k(A,B,C)\to B_R(A,B,C)
\end{equation}
\noindent
is a quasi-isomorphism of $R$CDGA's.   
\end{lemma}

\begin{proof}
We have that $B_R(A,B,C)=B_k(A,B,C)/V$, where $V$ is the sub-complex
generated by all elements of the form:

\begin{equation}
(a,b_1,\dots,rb_i,\dots, b_k, c)-(a, b_1,\dots, rb_{i+1}, \dots, b_k, c)
\end{equation}
where $r\in R, a\in A, b_j\in B$, and $c\in C$. It is due to the fact
that all elements of $R$ have even degree that $V$ is closed under
the differential $D=d+\delta$.  We claim that $V$ is, in fact, acyclic.
To see this, consider the map $s: V^i\to V^{i-1}$ defined by

\begin{align}
&s\{(a,b_1,\dots,rb_i,\dots, b_k,c)-(a,b_1,\dots,rb_{i+1},\dots, b_k,c)\}\\
&=(-1)^{\epsilon_i}\{(a,b_1,\dots,rb_i,1,b_{i+1},\dots, b_k,c)-
(a,b_1,\dots, b_i,r,b_{i+1},\dots,b_k,c)\}\notag
\end{align}

\noindent
where $\epsilon_i=deg\alpha+deg\omega_1+\dots +deg\omega_i-i$.  It is 
straightforward but tedious to check that $ds+sd=0$, and that
$\delta s+s\delta=id.$, so that $Ds+sD=id.$ and consequently $V$ is
acyclic.  Moreover, it is easy to check that $V$ is an ideal so that
the product on the bar complex induces a product on the bar complex over $R$.  

\end{proof}

\begin{remark}\label{lastremark}
Lemma \ref{overR} is valid using the normalized bar complex.
\end{remark}

\begin{corollary}
In the situation of \ref{eqlemma}, 
\begin{equation}
\theta: \bar B_R(A_G^\bu(X), A_G^\bu(B), A_G^\bu(E))\to A_G^\bu(E_f)
\end{equation}
\noindent
is a quasi-isomorphism of $R$CDGA$^o$'s ($R$CDGA's).  
\end{corollary}

\bibliographystyle{amsplain}
\bibliography{master}

\end{document}